\documentstyle[12pt]{article}
 \normalsize

\begin{document}
\title{{\bf {On the difference between a D. H. Lehmer number and its inverse over short interval}}
\footnotetext{\\
1.nyn0902@mail.nwpu.edu.cn \qquad 2.marong@nwpu.edu.cn \qquad 3.2805256907@qq.com}
\author{Yana Niu$^{1}$ \qquad Rong Ma$^{2}$ \qquad Haodong Wang$^{3}$
\\{\small{{School of Mathematics and Statistics, Northwestern
Polytechnical University}} }\\{\small{Xi'an, Shaanxi, 710072,
People's Republic of China}} \\}
\date{}}
\maketitle
\vspace{0.0cm}

\begin{center}
\large{{\bf Abstract}}
\end{center}

Let $q>2$ be an odd integer. For each integer $x$ with $0<x<q$ and $(q,x)= 1$, we know that there exists one and only one $\bar{x}$ with
$0<\bar{x}<q$ such that $x\bar{x}\equiv1(\bmod q)$. A Lehmer number is defined to be any integer $a$ with $2\dag(a+\bar{a})$. For any nonnegative integer $k$, Let
$$
M(x,q,k)=\displaystyle\mathop
{\displaystyle\mathop{\sum{'}}_{a=1}^{q}
\displaystyle\mathop{\sum{'}}_{b\leq xq}}
_{\mbox{\tiny$\begin{array}{c}
2|a+b+1\\
ab\equiv1(\bmod q)\end{array}$}}(a-b)^{2k}.
$$
The main purpose of this paper is to study the properties of $M(x,q,k)$, and give a sharp asymptotic formula, by using estimates of Kloosterman's sums and properties of trigonometric sums.

\vspace{1cm}
\begin{center}
\large{{\bf1. Introduction}}
\end{center}

Let $q>2$ be an odd integer. For each integer $x$ with $0<x<q$ and $(q,x)= 1$, we know that there exists one and only one $\bar{x}$ with
$0<\bar{x}<q$ such that $x\bar{x}\equiv1(\bmod q)$. Let $r(q)$ be the number of cases in which $x$ and $\bar{x}$ are of opposite parity. For $q=p$ a prime, D. H. Lehmer [1] asks us to find $r(p)$ or at least to say something nontrivial about it. About this problem, a lot of scholars [2,3] have studied it. For the sake of simplicity, we call such a number $x$ as a D. H. Lehmer number.

W. Zhang [4] has given an asymptotic estimate:
$$r(p)=\frac{1}{2}p+O(p^\frac{1}{2}\ln^2p).
$$
Later, W. Zhang [5,6] also proved that for every odd integer $q\geq3$,
$$r(q)=\frac{1}{2}\phi(q)+O(q^\frac{1}{2}d^2(q)\ln^2q),
$$
where $\phi(q)$ is the Euler function and $d(q)$ is the divisor function.

For any nonnegative integer $k$, let
$$
M(q,k)=\displaystyle\mathop{\displaystyle\mathop{\sum{'}}_{a=1}^{q}}_{2|a+\bar{a}+1}(a-\bar{a})^{2k},
$$
W. Zhang [7] gave a sharp asymptotic formula for $M(q,k)$ as following:
$$
M(q,k)=\frac{1}{(2k+1)(2k+2)}\phi(q)q^{2k}+O(4^{k}q^{2k+\frac{1}{2}}d^2(q)\ln^{2}q).
$$
where ${\sum{'}}_{a}$ denotes the summation over all $a$ such that $(a,q)=1$.
Moreover, he [8] also proved
$$
\displaystyle\mathop{\displaystyle\mathop{\sum}_{a=1}^{q}}_{(a,q)=1}(a-\bar{a})^{2k}=
\frac{1}{(2k+1)(k+1)}\phi(q)q^{2k}+O(4^{k}q^{2k+\frac{1}{2}}d^2(q)\ln^{2}q).
$$

The main purpose of this paper is to study the distribution properties of D. H. Lehmer numbers and the asymptotic properties of the 2kth power mean
$$
M(x,q,k)=\displaystyle\mathop
{\displaystyle\mathop{\sum{'}}_{a=1}^{q}
\displaystyle\mathop{\sum{'}}_{b\leq xq}}
_{\mbox{\tiny$\begin{array}{c}
2|a+b+1\\
ab\equiv1(\bmod q)\end{array}$}}(a-b)^{2k}.
$$
It seems that no one has studied this problem yet. The problem is interesting because it can help us to find how large is the difference between a D. H. Lehmer number and its inverse modulo $q$. In this paper, we use estimates of Kloosterman¡¯s sums and properties of trigonometric sums to give a sharper asymptotic formula for $M(x,q,k)$ for any fixed positive integer $k$. That is, we shall prove the following:

\noindent \textbf{Theorem. }For any odd number $q$ and integer $k$, we have the asymptotic formula
$$
M(x,q,k)=\frac{1}{(2k+1)(2k+2)}x\phi(q)q^{2k}+O(q^{2k+\frac{1}{2}}\ln^{2}q)
$$
where $\phi(q)$ is the Euler function.

\begin{center}
\large{{\bf 2. Some lemmas}}
\end{center}

In this section, we prove some elementary lemmas which are necessary in the proof of the theorems.

\noindent \textbf{Lemma 1. }Let $q$ be an odd number. For any integer $n$ and nonnegative integer $r$, define
$$
K(n,r)=\sum_{a=1}^{q}a^{r}e\left(\frac{an}{q}\right), \qquad H(n,r)=\sum_{a=1}^{q}(-1)^{a}a^{r}e\left(\frac{an}{q}\right),
$$
where $e(y)=e^{2\pi iy}$. We have the estimates
\begin{eqnarray}
K(n,r)\,\,\,\,\left\{
\begin{array}{ll}
= \frac{q^{r+1}}{r+1}+O(q^r),& q|n\\
\ll \frac{q^r}{|\sin(\pi n/q)|},& q\dag n
\end{array}
\right.
\end{eqnarray}
\begin{eqnarray}
H(n,r)\ll \frac{q^r}{|\cos(\pi n/q)|}.
\end{eqnarray}
\\
\textbf{Proof. }See Ref. [7].

\noindent \textbf{Lemma 2. }For any integer $K\geq1$ and $0<\alpha<1$, we have
$$
\left|\sum_{n=1}^{K}e(\alpha n)\right|\leq \min\left(K,\frac{1}{2\langle\alpha\rangle}\right),
$$
where $\langle\alpha\rangle=\min(\{\alpha\},1-\{\alpha\})$, $\{\alpha\}$ is the decimal part of $\alpha$. \\
\\
\textbf{Proof. }See Ref. [9].

\noindent \textbf{Lemma 3. }For any integer $q\geq3$, $x>\frac{1}{2}$, and any positive integer $n\geq1$, we have
$$
\sum_{l=1}^{q-1}\left|\displaystyle\mathop{\sum{'}}_{b\leq xq}e\left(\frac{b(n-l)}{q}\right)\right|\ll q^{1+\epsilon}.
$$
\\
\textbf{Proof. }From Lemma 2, when $n\not\equiv0(\bmod q)$, for $1\leq l\leq q-1$, there must be one and only one $l$ such that $n-l\equiv0(\bmod q)$, so we get
\begin{eqnarray}
&&\sum_{l=1}^{q-1}\left|\displaystyle\mathop{\sum{'}}_{b\leq xq}e\left(\frac{b(n-l)}{q}\right)\right|\nonumber\\
&=&\sum_{l=1}^{q-1}\left|\sum_{b\leq xq}e\left(\frac{b(n-l)}{q}\right)\sum_{d|(b,q)}\mu(d)\right|\nonumber\\
&=&\sum_{l=1}^{q-1}\left|\sum_{d|q}\mu(d)\displaystyle\mathop{\displaystyle\mathop{\sum}_{b\leq xq}}_{d|b}e\left(\frac{b(n-l)}{q}\right)\right|\nonumber\\
&=&\displaystyle\mathop
{\displaystyle\mathop{\sum}_{l=1}^{q-1}}_{l\not\equiv n(\bmod q)}\left|\sum_{d|q}\mu(d)\displaystyle\mathop{\sum}_{b\leq xq/d}e\left(\frac{b(n-l)}{q}\right)\right|+\left|\sum_{d|q}\mu(d)\displaystyle\mathop{\sum}_{b\leq xq/d}1\right|\nonumber\\
&\leq&\displaystyle\mathop
{\displaystyle\mathop{\sum}_{l=1}^{q-1}}_{l\not\equiv n(\bmod q)}\sum_{d|q}\left|\min\left\{\frac{xq}{d},\frac{1}{2\langle\frac{n-l}{q}\rangle}\right\}\right|
+\left|xq\sum_{d|q}\frac{\mu(d)}{d}\right|\nonumber\\
&\leq&\displaystyle\mathop
{\displaystyle\mathop{\sum}_{l=1}^{q-1}}_{l\not\equiv n(\bmod q)}\sum_{d|q}\left|\frac{1}{2\langle\frac{n-l}{q}\rangle}\right|+x\phi(q)\nonumber\\
&=&\sum_{l=1}^{q-1}\sum_{d|q}\left|\frac{1}{2\langle\frac{l}{q}\rangle}\right|-
\sum_{d|q}\left|\frac{1}{2\langle\frac{n}{q}\rangle}\right|+x\phi(q)\nonumber\\
&\leq&2\sum_{1\leq l\leq[\frac{q}{2}]}\sum_{d|q}\frac{1}{\frac{2l}{q}}-
\sum_{d|q}\left|\frac{1}{2\langle\frac{n}{q}\rangle}\right|+x\phi(q)\nonumber\\
&=&q\sum_{1\leq l\leq[\frac{q}{2}]}\sum_{d|q}\frac{1}{l}-
\sum_{d|q}\left|\frac{1}{2\langle\frac{n}{q}\rangle}\right|+x\phi(q)\nonumber\\
&=&qd(q)\sum_{1\leq l\leq[\frac{q}{2}]}\frac{1}{l}-d(q)\left|\frac{1}{2\langle\frac{n}{q}\rangle}\right|+x\phi(q)\nonumber\\
&\ll&q^{1+\epsilon}.
\end{eqnarray}

When $n\equiv0(\bmod q)$, from Lemma 2, we also get
\begin{eqnarray}
&&\sum_{l=1}^{q-1}\left|\displaystyle\mathop{\sum{'}}_{b\leq xq}e\left(\frac{b(n-l)}{q}\right)\right|\nonumber\\
&=&\sum_{l=1}^{q-1}\left|\sum_{b\leq xq}e\left(\frac{b(n-l)}{q}\right)\sum_{d|(b,q)}\mu(d)\right|\nonumber\\
&=&\sum_{l=1}^{q-1}\left|\sum_{d|q}\mu(d)\displaystyle\mathop{\sum}_{b\leq xq/d}e\left(\frac{bl}{q}\right)\right|\nonumber\\
&\leq&\sum_{l=1}^{q-1}\sum_{d|q} \left|\min\left\{\frac{xq}{d},\frac{1}{2\langle\frac{l}{q}\rangle}\right\}\right|\nonumber\\
&\leq&\sum_{l=1}^{q-1}\sum_{d|q}\left|\frac{1}{2\langle\frac{l}{q}\rangle}\right|\nonumber\\
&\leq&2\sum_{1\leq l\leq[\frac{q}{2}]}\sum_{d|q}\frac{1}{\frac{2l}{q}}\nonumber\\
&=&q\sum_{1\leq l\leq[\frac{q}{2}]}\sum_{d|q}\frac{1}{l}\nonumber\\
&=&qd(q)\sum_{1\leq l\leq[\frac{q}{2}]}\frac{1}{l}\nonumber\\
&\ll&q^{1+\epsilon}.
\end{eqnarray}

Therefore, combining (3) and (4), for any positive integer $n\geq1$, we have the estimate
$$
\sum_{l=1}^{q-1}\left|\displaystyle\mathop{\sum{'}}_{b\leq xq}e\left(\frac{b(n-l)}{q}\right)\right|\ll q^{1+\epsilon}.
$$
This proves Lemma 3.

\noindent \textbf{Lemma 4. }Let $q\geq3$ be an integer, and $\chi$ denote the Dirichlet character modulo $q$. The Gauss sum is defined by
$$G(n,\chi)=\sum_{a=1}^{q}\chi(a)e\left(\frac{na}{q}\right).
$$
Then by the principal Dirichlet character $\chi_0$ modulo $q$, we have the identity
$$G(n,\chi_0)=\mu\left(\frac{q}{(n,q)}\right)\phi(q)\phi^{-1}\left(\frac{q}{(n,q)}\right),
$$
where $\mu(n)$ is the M\"{o}bius function, $\phi(q)$ is the Euler function, and $(n,q)$ is the greatest common divisor of $n$ and $q$.\\
\\
\textbf{Proof. }See Ref. [10].

\noindent \textbf{Lemma 5. }For $q>2$ an integer, any non-principal Dirichlet character $\chi$ modulo $q$ and any positive integers $m$ and $l$, we have the estimate
\begin{eqnarray}
\sum_{\chi\bmod q}G(m,\chi)G(l,\chi)\ll(m,l,q)^{\frac{1}{2}}q^{\frac{3}{2}+\epsilon},
\end{eqnarray}
where $(m,l,q)$ is the greatest common divisor of $m$, $l$ and $q$, $\epsilon$ is any fixed positive real number.

On the other hand, for the principal Dirichlet character $\chi_0$, any positive integers $m$ and $l$, we also have
\begin{eqnarray}
G(m,\chi_0)G(l,\chi_0)\ll(m,q)(l,q)q^{\epsilon}.
\end{eqnarray}
\\
\textbf{Proof. }First, we prove (5). According to the orthogonality of character sums, we have
\begin{eqnarray}
\sum_{\chi\bmod q}\chi(n)=\left\{
\begin{array}{ll}
\phi(q),& n\equiv l(\bmod q)\\
0,& n\not\equiv l(\bmod q)
\end{array}
\right.\nonumber
\end{eqnarray}
and hence we have
\begin{eqnarray}
&&\sum_{\chi\bmod q}G(m,\chi)G(l,\chi)\nonumber\\
&=&\sum_{\chi\bmod q}\sum_{s=1}^{q}\chi(s)e\left(\frac{ms}{q}\right)
\sum_{t=1}^{q}\chi(t)e\left(\frac{lt}{q}\right)\nonumber\\
&=&\displaystyle\mathop{\displaystyle\mathop{\sum{'}}_{s=1}^{q}
\displaystyle\mathop{\sum{'}}_{t=1}^{q}}e\left(\frac{ms+lt}{q}\right)\sum_{\chi\bmod q}\chi(s)\chi(t)\nonumber\\
&=&\phi(q)\displaystyle\mathop
{\displaystyle\mathop{\sum}_{s=1}^{q-1}
\displaystyle\mathop{\sum}_{t=1}^{q-1}}
_{st\equiv1(\bmod q)}e\left(\frac{ms+lt}{q}\right)\nonumber\\
&=&\phi(q)\displaystyle\mathop{\sum{'}}_{t=1}^{q-1}e\left(\frac{m\bar{t}+lt}{q}\right)\nonumber\\
&\ll&\phi(q)(m,l,q)^{\frac{1}{2}}q^{\frac{1}{2}}\nonumber\\
&\ll&(m,l,q)^{\frac{1}{2}}q^{\frac{3}{2}+\epsilon}.\nonumber
\end{eqnarray}

Now we show (6). From Lemma 4 and the definition of Gauss sum, we have
\begin{eqnarray}
&&G(m,\chi_0)G(l,\chi_0)\nonumber\\
&=&\mu\left(\frac{q}{(m,q)}\right)\phi(q)\phi^{-1}\left(\frac{q}{(m,q)}\right)
\mu\left(\frac{q}{(l,q)}\right)\phi(q)\phi^{-1}\left(\frac{q}{(l,q)}\right)\nonumber\\
&\ll&\phi^2(q)\frac{(m,q)(l,q)d^2(q)}{q^2}\nonumber\\
&\ll&(m,q)(l,q)q^{\epsilon},\nonumber
\end{eqnarray}
where we have used $\phi(q)\gg\frac{q}{d(q)}$(see Ref. [11]), and $d(q)$ is the divisor function.

\noindent \textbf{Lemma 6. }Let $m$, $n$ and $q$ be integers, and $q>2$. Then we have the estimates
$$
S(x,m,n;q)=\displaystyle\mathop
{\displaystyle\mathop{\sum{'}}_{a=1}^{q}
\displaystyle\mathop{\sum{'}}_{b\leq xq}}
_{ab\equiv1(\bmod q)}e\left(\frac{am+bn}{q}\right)\ll q^{\frac{1}{2}+\epsilon}
$$
\\
\textbf{Proof. }From the orthogonality of character sums, we have
\begin{eqnarray}
&&S(x,m,n;q)\nonumber\\
&=&\displaystyle\mathop
{\displaystyle\mathop{\sum{'}}_{a=1}^{q}
\displaystyle\mathop{\sum{'}}_{b\leq xq}}
_{ab\equiv1(\bmod q)}e\left(\frac{am+bn}{q}\right)\nonumber\\
&=&\frac{1}{\phi(q)}\displaystyle\mathop
{\displaystyle\mathop{\sum{'}}_{a=1}^{q}\displaystyle\mathop{\sum{'}}_{b\leq xq}}e\left(\frac{am+bn}{q}\right)\sum_{\chi\bmod q}\chi(a)\chi(b)\nonumber\\
&=&\frac{1}{\phi(q)}\sum_{\chi\bmod q}\left(\displaystyle\mathop{\sum{'}}_{a=1}^{q}\chi(a)e\left(\frac{am}{q}\right)\right)
\left(\displaystyle\mathop{\sum{'}}_{b\leq xq}\chi(b)e\left(\frac{bn}{q}\right)\right)\nonumber\\
&=&\frac{1}{\phi(q)}\sum_{\chi\neq\chi_0}\left(\displaystyle\mathop{\sum{'}}_{a=1}^{q}\chi(a)e\left(\frac{am}{q}\right)\right)
\left(\displaystyle\mathop{\sum{'}}_{b\leq xq}\chi(b)e\left(\frac{bn}{q}\right)\right)\nonumber\\
&&+\frac{1}{\phi(q)}\left(\displaystyle\mathop{\sum{'}}_{a=1}^{q-1}e\left(\frac{am}{q}\right)\right)
\left(\displaystyle\mathop{\sum{'}}_{b\leq xq}e\left(\frac{bn}{q}\right)\right)\nonumber\\
&=&S_1+S_2.
\end{eqnarray}
Now we will estimate both $S_1$ and $S_2$ respectively. Firstly, we shall estimate $S_1$. From the identity, for any Dirichlet character $\chi\neq\chi_0$ modulo $q$,
$$\chi(a)=\frac{1}{q}\sum_{k=1}^{q-1}G(k,\chi)e\left(-\frac{ak}{q}\right).
$$
Hence according to Lemma 3 and Lemma 5, we have
\begin{eqnarray}
S_1&=&\frac{1}{\phi(q)}\sum_{\chi\neq\chi_0}
\left(\displaystyle\mathop{\sum{'}}_{a=1}^{q}\frac{1}{q}\sum_{k=1}^{q-1}G(k,\chi)e\left(-\frac{ak}{q}\right)e\left(\frac{am}{q}\right)\right)\times\nonumber\\
&&\times\left(\displaystyle\mathop{\sum{'}}_{b\leq xq}\frac{1}{q}\sum_{l=1}^{q-1}G(l,\chi)e\left(-\frac{bl}{q}\right)e\left(\frac{bn}{q}\right)\right)\nonumber\\
&=&\frac{1}{q^2\phi(q)}\sum_{\chi\neq\chi_0}
\left(\sum_{k=1}^{q-1}G(k,\chi)\displaystyle\mathop{\sum{'}}_{a=1}^{q}e\left(\frac{a(m-k)}{q}\right)\right)\times\nonumber\\
&&\times\left(\sum_{l=1}^{q-1}G(l,\chi)\displaystyle\mathop{\sum{'}}_{b\leq xq}e\left(\frac{b(n-l)}{q}\right)\right)\nonumber\\
&=&\frac{1}{q^2\phi(q)}\sum_{\chi\neq\chi_0}
\left(\sum_{k=1}^{q-1}G(k,\chi)\sum_{a=1}^{q}e\left(\frac{a(m-k)}{q}\right)\sum_{d|(a,q)}\mu(d)\right)\times\nonumber\\
&&\times\left(\sum_{l=1}^{q-1}G(l,\chi)\displaystyle\mathop{\sum{'}}_{b\leq xq}e\left(\frac{b(n-l)}{q}\right)\right)\nonumber\\
&=&\frac{1}{q^2\phi(q)}\sum_{\chi\neq\chi_0}
\left(\sum_{k=1}^{q-1}G(k,\chi)\sum_{d|q}\mu(d)\displaystyle\mathop{\displaystyle\mathop{\sum}_{a=1}^{q}}_{d|a}e\left(\frac{a(m-k)}{q}\right)\right)\times\nonumber\\
&&\times\left(\sum_{l=1}^{q-1}G(l,\chi)\displaystyle\mathop{\sum{'}}_{b\leq xq}e\left(\frac{b(n-l)}{q}\right)\right)\nonumber\\
&=&\frac{q}{q^2\phi(q)}\sum_{\chi\neq\chi_0}
\left(\displaystyle\mathop{\displaystyle\mathop{\sum}_{k=1}^{q-1}}_{q|(m-k)}G(k,\chi)\sum_{d|q}\frac{\mu(d)}{d}\right)\times\nonumber\\
&&\times\left(\sum_{l=1}^{q-1}G(l,\chi)\displaystyle\mathop{\sum{'}}_{b\leq xq}e\left(\frac{b(n-l)}{q}\right)\right)\nonumber\\
&=&\frac{1}{q^2}\sum_{\chi\neq\chi_0}G(m,\chi)
\left(\sum_{l=1}^{q-1}G(l,\chi)\displaystyle\mathop{\sum{'}}_{b\leq xq}e\left(\frac{b(n-l)}{q}\right)\right)\nonumber\\
&=&\frac{1}{q^2}
\sum_{l=1}^{q-1}\displaystyle\mathop{\sum{'}}_{b\leq xq}e\left(\frac{b(n-l)}{q}\right)\sum_{\chi\neq\chi_0}G(m,\chi)G(l,\chi)\nonumber\\
&=&\frac{1}{q^2}\sum_{l=1}^{q-1}\displaystyle\mathop{\sum{'}}_{b\leq xq}e\left(\frac{b(n-l)}{q}\right)\sum_{\chi\bmod q}G(m,\chi)G(l,\chi)+\nonumber\\
&&+\frac{1}{q^2}\sum_{l=1}^{q-1}\displaystyle\mathop{\sum{'}}_{b\leq xq}e\left(\frac{b(n-l)}{q}\right)G(m,\chi_0)G(l,\chi_0)\nonumber\\
&\leq&\frac{1}{q^2}\sum_{l=1}^{q-1}\left|\displaystyle\mathop{\sum{'}}_{b\leq xq}e\left(\frac{b(n-l)}{q}\right)\right|\left|\sum_{\chi\bmod q}G(m,\chi)G(l,\chi)\right|+\nonumber\\
&&+\frac{1}{q^2}\sum_{l=1}^{q-1}\left|\displaystyle\mathop{\sum{'}}_{b\leq xq}e\left(\frac{b(n-l)}{q}\right)\right||G(m,\chi_0)G(l,\chi_0)|\nonumber\\
&\leq&\frac{1}{q^2}\sum_{l=1}^{q-1}\left|\displaystyle\mathop{\sum{'}}_{b\leq xq}e\left(\frac{b(n-l)}{q}\right)\right|(m,l,q)^{\frac{1}{2}}q^{\frac{3}{2}+\epsilon}+\nonumber\\
&&+\frac{1}{q\phi(q)}\sum_{l=1}^{q-1}\left|\displaystyle\mathop{\sum{'}}_{b\leq xq}e\left(\frac{b(n-l)}{q}\right)\right|(m,q)(l,q)q^{\epsilon}\nonumber\\
&\ll&q^{-\frac{1}{2}+\epsilon}\sum_{l=1}^{q-1}\left|\displaystyle\mathop{\sum{'}}_{b\leq xq}e\left(\frac{b(n-l)}{q}\right)\right|\nonumber\\
&\ll&q^{\frac{1}{2}+\epsilon}.
\end{eqnarray}

Now we estimate $S_2$, we have
\begin{eqnarray}
|S_2|&=&\left|\frac{1}{\phi(q)}\left(\displaystyle\mathop{\sum{'}}_{a=1}^{q-1}e\left(\frac{am}{q}\right)\right)
\left(\displaystyle\mathop{\sum{'}}_{b\leq xq}e\left(\frac{bn}{q}\right)\right)\right|\nonumber\\
&=&\frac{1}{\phi(q)}\left|\displaystyle\mathop{\sum{'}}_{a=1}^{q-1}e\left(\frac{am}{q}\right)\right|
\left|\displaystyle\mathop{\sum{'}}_{b\leq xq}e\left(\frac{bn}{q}\right)\right|\nonumber\\
&\leq&\frac{1}{\phi(q)}\left|\displaystyle\mathop{\sum{'}}_{b\leq xq}e\left(\frac{bn}{q}\right)\right|\nonumber\\
&\leq&\frac{xq}{\phi(q)}\nonumber\\
&\leq&xq^{\epsilon}\nonumber\\
&\ll&q^{\frac{1}{2}+\epsilon},
\end{eqnarray}
where $x<q^{\frac{1}{2}+\epsilon}$.

Therefore, from (7)-(9), we have
$$
S(x,m,n;q)=\displaystyle\mathop
{\displaystyle\mathop{\sum{'}}_{a=1}^{q}
\displaystyle\mathop{\sum{'}}_{b\leq xq}}
_{ab\equiv1(\bmod q)}e\left(\frac{am+bn}{q}\right)
\ll q^{\frac{1}{2}+\epsilon},
$$
where $\epsilon$ is any positive real number.

\noindent \textbf{Lemma 7. }Let $r$, $s$ and $q$ be positive integers and $q>2$. Then
$$
\displaystyle\mathop
{\displaystyle\mathop{\sum{'}}_{a=1}^{q}
\displaystyle\mathop{\sum{'}}_{b\leq xq}}
_{ab\equiv1(\bmod q)}a^{r}b^{s}=\frac{x\phi(q)q^{r+s}}{(r+1)(s+1)}+O(q^{r+s+\frac{1}{2}}\ln^2q),
$$
where $\phi(q)$ is the Euler function.\\
\\
\textbf{Proof. }From the trigonometric identity
$$
\sum_{a=1}^{q}e\left(\frac{an}{q}\right)=\left\{
\begin{array}{ll}
q,& q|n\\
0,& q\dag n
\end{array}
\right.
$$
we get the identity
\begin{eqnarray}
&&\displaystyle\mathop
{\displaystyle\mathop{\sum{'}}_{a=1}^{q}
\displaystyle\mathop{\sum{'}}_{b\leq xq}}
_{ab\equiv1(\bmod q)}a^{r}b^{s}\nonumber\\
&=&\frac{1}{q^2}\displaystyle\mathop
{\displaystyle\mathop{\sum{'}}_{a=1}^{q}
\displaystyle\mathop{\sum{'}}_{b\leq xq}}_{ab\equiv1(\bmod q)}\sum_{c,d=1}^{q}c^{r}d^{s}\sum_{m,n=1}^{q}e\left(\frac{m(a-c)+n(b-d)}{q}\right)\nonumber\\
&=&\frac{1}{q^2}\sum_{m,n=1}^{q}\left(\displaystyle\mathop
{\displaystyle\mathop{\sum{'}}_{a=1}^{q}
\displaystyle\mathop{\sum{'}}_{b\leq xq}}
_{ab\equiv1(\bmod q)}e\left(\frac{am+bn}{q}\right)\right)
\left(\sum_{c=1}^{q}c^{r}e\left(\frac{-mc}{q}\right)\right)
\left(\sum_{d=1}^{q}d^{s}e\left(\frac{-nd}{q}\right)\right)\nonumber\\
&=&\frac{1}{q^2}\sum_{m=1}^{q}\sum_{n=1}^{q}S(x,m,n;q)K(-m,r)K(-n,s)\nonumber\\
&=&\frac{1}{q^2}\sum_{m=1}^{q-1}S(x,m,q;q)K(-m,r)K(-q,s)+\nonumber\\
&&+\frac{1}{q^2}\sum_{n=1}^{q-1}S(x,q,n;q)K(-q,r)K(-n,s)+\nonumber\\
&&+\frac{1}{q^2}\sum_{m=1}^{q-1}\sum_{n=1}^{q-1}S(x,m,n;q)K(-m,r)K(-n,s)+\nonumber\\
&&+\frac{1}{q^2}S(x,q,q;q)K(-q,r)K(-q,s),
\end{eqnarray}
where $K(-m,r)$ is defined in Lemma 1. From (2) of Lemma 1, Lemma 6 and noting that $2/\pi\leq(\sin x/x)$ for $|x|\leq \pi/2$, we get
\begin{eqnarray}
&&\frac{1}{q^2}S(x,q,q;q)K(-q,r)K(-q,s)\nonumber\\
&=&\frac{1}{q^2}\left(\displaystyle\mathop
{\displaystyle\mathop{\sum{'}}_{a=1}^{q}
\displaystyle\mathop{\sum{'}}_{b\leq xq}}
_{ab\equiv1(\bmod q)}e\left(\frac{aq+bq}{q}\right)\right)
\left(\frac{q^{r+1}}{r+1}+O(q^r)\right)\left(\frac{q^{s+1}}{s+1}+O(q^s)\right)\nonumber\\
&=&\frac{1}{q^2}\left(\displaystyle\mathop
{\displaystyle\mathop{\sum{'}}_{a=1}^{q}
\displaystyle\mathop{\sum{'}}_{b\leq xq}}
_{ab\equiv1(\bmod q)}1\right)
\left(\frac{q^{r+1}}{r+1}+O(q^r)\right)\left(\frac{q^{s+1}}{s+1}+O(q^s)\right)\nonumber\\
&=&\frac{1}{q^2}\left(\displaystyle\mathop{\sum{'}}_{b\leq xq}1\right)
\left(\frac{q^{r+1}}{r+1}+O(q^r)\right)\left(\frac{q^{s+1}}{s+1}+O(q^s)\right)\nonumber\\
&=&\frac{1}{q^2}\left(\sum_{b\leq xq}\sum_{d|(b,q)}\mu(d)\right)
\left(\frac{q^{r+1}}{r+1}+O(q^r)\right)\left(\frac{q^{s+1}}{s+1}+O(q^s)\right)\nonumber\\
&=&\frac{1}{q^2}\left(\sum_{d|q}\mu(d)\sum_{b\leq xq/d}1\right)
\left(\frac{q^{r+1}}{r+1}+O(q^r)\right)\left(\frac{q^{s+1}}{s+1}+O(q^s)\right)\nonumber\\
&=&\frac{1}{q^2}\left(xq\sum_{d|q}\frac{\mu(d)}{d}\right)
\left(\frac{q^{r+1}}{r+1}+O(q^r)\right)\left(\frac{q^{s+1}}{s+1}+O(q^s)\right)\nonumber\\
&=&\frac{x\phi(q)}{q^2}\left(\frac{q^{r+1}}{r+1}+O(q^r)\right)\left(\frac{q^{s+1}}{s+1}+O(q^s)\right)\nonumber\\
&=&\frac{x\phi(q)q^{r+s}}{(r+1)(s+1)}+O(q^{r+s}),
\end{eqnarray}
\begin{eqnarray}
&&\sum_{m=1}^{q-1}S(x,m,q;q)K(-m,r)K(-q,s)\nonumber\\
&\ll&\sum_{m=1}^{q-1}q^{\frac{1}{2}}q^{s+1}\frac{q^r}{|\sin(\pi m/q)|}\nonumber\\
&\ll&\sum_{m=1}^{q-1}q^{\frac{1}{2}}q^{s+1}q^r\frac{q}{2m}\nonumber\\
&\ll&q^{r+s+\frac{5}{2}}\sum_{m=1}^{q-1}\frac{1}{m}\nonumber\\
&\ll&q^{r+s+\frac{5}{2}}\ln q.
\end{eqnarray}
Similarly, we can get the estimate
\begin{eqnarray}
&&\sum_{n=1}^{q-1}S(x,q,n;q)K(-q,r)K(-n,s)\ll q^{r+s+\frac{5}{2}}\ln q.
\end{eqnarray}
\begin{eqnarray}
&&\sum_{m=1}^{q-1}\sum_{n=1}^{q-1}S(x,m,n;q)K(-m,r)K(-n,s)\nonumber\\
&\ll&\sum_{m=1}^{q-1}\sum_{n=1}^{q-1}\frac{q^{\frac{1}{2}}q^{r+s}}{|\sin(\pi m/q)||\sin(\pi n/q)|}\nonumber\\
&\ll&q^{r+s+\frac{5}{2}}\sum_{m=1}^{q-1}\sum_{n=1}^{q-1}\frac{1}{mn}\nonumber\\
&\ll&q^{r+s+\frac{5}{2}}\ln^2q.
\end{eqnarray}

Combining (10)-(14) we immediately deduce that
$$
\displaystyle\mathop
{\displaystyle\mathop{\sum{'}}_{a=1}^{q}
\displaystyle\mathop{\sum{'}}_{b\leq xq}}
_{ab\equiv1(\bmod q)}a^{r}b^{s}=\frac{x\phi(q)q^{r+s}}{(r+1)(s+1)}+O(q^{r+s+\frac{1}{2}}\ln^2q).
$$
This is the conclusion of Lemma 7.

\noindent \textbf{Lemma 8. }Let $r$, $s$ and $q$ be positive integers and $q>2$. Then
$$
\displaystyle\mathop
{\displaystyle\mathop{\sum{'}}_{a=1}^{q}
\displaystyle\mathop{\sum{'}}_{b\leq xq}}
_{ab\equiv1(\bmod q)}(-1)^{a+b}a^{r}b^{s}=O(q^{r+s+\frac{1}{2}}\ln^2q).
$$
\\
\textbf{Proof. }Similarly, we get
\begin{eqnarray}
&&\displaystyle\mathop
{\displaystyle\mathop{\sum{'}}_{a=1}^{q}
\displaystyle\mathop{\sum{'}}_{b\leq xq}}
_{ab\equiv1(\bmod q)}(-1)^{a+b}a^{r}b^{s}\nonumber\\
&=&\frac{1}{q^2}\displaystyle\mathop
{\displaystyle\mathop{\sum{'}}_{a=1}^{q}
\displaystyle\mathop{\sum{'}}_{b\leq xq}}
_{ab\equiv1(\bmod q)}\sum_{c,d=1}^{q}(-1)^{c+d}c^{r}d^{s}\sum_{m,n=1}^{q}e\left(\frac{m(a-c)+n(b-d)}{q}\right)\nonumber\\
&=&\frac{1}{q^2}\sum_{m,n=1}^{q}\left(\displaystyle\mathop
{\displaystyle\mathop{\sum{'}}_{a=1}^{q}
\displaystyle\mathop{\sum{'}}_{b\leq xq}}
_{ab\equiv1(\bmod q)}e\left(\frac{am+bn}{q}\right)\right)
\left(\sum_{c=1}^{q}(-1)^{c}c^{r}e\left(\frac{-mc}{q}\right)\right)
\left(\sum_{d=1}^{q}(-1)^{d}d^{s}e\left(\frac{-nd}{q}\right)\right)\nonumber\\
&=&\frac{1}{q^2}\sum_{m=1}^{q}\sum_{n=1}^{q}S(x,m,n;q)H(-m,r)H(-n,s)\nonumber\\
&\ll&\frac{1}{q^2}\sum_{m=1}^{q}\sum_{n=1}^{q}\frac{q^{\frac{1}{2}}q^{r+s}}
{|\cos(\pi m/q)||\cos(\pi n/q)|}\nonumber\\
&\ll&\frac{1}{q^2}\sum_{m=1}^{q}\sum_{n=1}^{q}\frac{q^{\frac{1}{2}}q^{r+s}q^2}
{(q-2m)(q-2n)}\nonumber\\
&=&q^{r+s+\frac{1}{2}}\sum_{m=1}^{q}\sum_{n=1}^{q}\frac{1}{(q-2m)(q-2n)}\nonumber\\
&\ll&q^{r+s+\frac{1}{2}}\ln^2q.\nonumber
\end{eqnarray}
where $H(-m,r)$ is defined in Lemma 1.

Noting that $|\cos(\pi m/q)|=|\sin(\pi(q-2m)/(2q))|$ and $q-2m\neq0$, and we have used (2) of Lemma 1 in the proof above.

\begin{center}
\large{{\bf 3. Proof of The Theorem}}
\end{center}

In this section, we shall complete the proof of the theorem. By the binomial
formula, Lemma 7 and Lemma 8 we get
\begin{eqnarray}
&&M(x,q,k)=\displaystyle\mathop
{\displaystyle\mathop{\sum{'}}_{a=1}^{q}
\displaystyle\mathop{\sum{'}}_{b\leq xq}}
_{\mbox{\tiny$\begin{array}{c}
2|a+b+1\\
ab\equiv1(\bmod q)\end{array}$}}(a-b)^{2k}\nonumber\\
&=&\frac{1}{2}\displaystyle\mathop
{\displaystyle\mathop{\sum{'}}_{a=1}^{q}
\displaystyle\mathop{\sum{'}}_{b\leq xq}}
_{ab\equiv1(\bmod q)}(1-(-1)^{a+b})(a-b)^{2k}\nonumber\\
&=&\frac{1}{2}\displaystyle\mathop
{\displaystyle\mathop{\sum{'}}_{a=1}^{q}
\displaystyle\mathop{\sum{'}}_{b\leq xq}}
_{ab\equiv1(\bmod q)}(a-b)^{2k}-
\frac{1}{2}\displaystyle\mathop
{\displaystyle\mathop{\sum{'}}_{a=1}^{q}
\displaystyle\mathop{\sum{'}}_{b\leq xq}}
_{ab\equiv1(\bmod q)}(-1)^{a+b}(a-b)^{2k}\nonumber\\
&=&\frac{1}{2}\sum_{i=0}^{2k}C_{2k}^{i}(-1)^i
\left(\displaystyle\mathop
{\displaystyle\mathop{\sum{'}}_{a=1}^{q}
\displaystyle\mathop{\sum{'}}_{b\leq xq}}
_{ab\equiv1(\bmod q)}a^{2k-i}b^{i}-
\displaystyle\mathop
{\displaystyle\mathop{\sum{'}}_{a=1}^{q}
\displaystyle\mathop{\sum{'}}_{b\leq xq}}
_{ab\equiv1(\bmod q)}(-1)^{a+b}a^{2k-i}b^{i}\right)\nonumber\\
&=&\frac{1}{2}\sum_{i=0}^{2k}C_{2k}^{i}(-1)^i
\left(\frac{x\phi(q)q^{2k}}{(i+1)(2k-i+1)}+O(q^{2k+\frac{1}{2}}\ln^2q)\right)+\nonumber\\
&&+O\left(\sum_{i=0}^{2k}C_{2k}^{i}q^{2k+\frac{1}{2}}\ln^2q\right)\nonumber\\
&=&\frac{x\phi(q)q^{2k}}{2}\sum_{i=0}^{2k}\frac{C_{2k}^{i}(-1)^i}{(i+1)(2k-i+1)}+
O\left(q^{2k+\frac{1}{2}}\ln^2q\right)\nonumber\\
&=&\frac{x\phi(q)q^{2k}}{2(2k+1)(2k+2)}\sum_{i=0}^{2k}(-1)^iC_{2k+2}^{i+1}+
O\left(q^{2k+\frac{1}{2}}\ln^2q\right)\nonumber\\
&=&\frac{x\phi(q)q^{2k}}{2(2k+1)(2k+2)}\left(-\sum_{i=0}^{2k+2}(-1)^iC_{2k+2}^{i}+2\right)+
O\left(q^{2k+\frac{1}{2}}\ln^2q\right)\nonumber\\
&=&\frac{x\phi(q)q^{2k}}{2(2k+1)(2k+2)}\left(-(1-1)^{2k+2}+2\right)+
O\left(q^{2k+\frac{1}{2}}\ln^2q\right)\nonumber\\
&=&\frac{x\phi(q)q^{2k}}{(2k+1)(2k+2)}+
O\left(q^{2k+\frac{1}{2}}\ln^2q\right).\nonumber
\end{eqnarray}
This is the conclusion of Theorem.
\\

\end{document}